\documentclass[leqno]{svjour3}
\usepackage{amsfonts}
\usepackage{amsmath}
\usepackage[colorlinks=true, pdfstartview=FitV, linkcolor=blue,citecolor=blue,
urlcolor=blue]{hyperref}
\usepackage{color}
\usepackage{amssymb}
\usepackage{graphicx}
\usepackage{ntheorem}

\smartqed
\theoremstyle{definition}
\newtheorem{lem}[theorem]{L\ e\ m\ m\ a}
\newtheorem{cor}[theorem]{C\ o\ r\ o\ l\ l\ a\ r\ y}
\newtheorem{prop}[theorem]{P\ r\ o\ p\ o\ s\ i\ t\ i\ o\ n}
\newtheorem{rem}[theorem]{R\ e\ m\ a\ r\ k}

\begin{document}

\title{Geometric Polynomials: Properties and Applications to Series with
Zeta Values}
\author{Khristo N. Boyadzhiev \and Ayhan Dil}
\institute{Khristo N. Boyadzhiev \at
              Department of Mathematics and Statistics, Ohio Northern University Ada, Ohio 45810, USA\\
              \email{k-boyadzhiev@onu.edu},
              \and
              Ayhan Dil \at
              Department of Mathematics, Akdeniz University, 07058-Antalya, Turkey\\
              \email{adil@akdeniz.edu.tr}
              }
\date{}
\maketitle

\begin{abstract}
We provide several properties of the geometric polynomials discussed in
earlier works of the authors. Further, the geometric polynomials are used to
obtain a closed form evaluation of certain series involving Riemann's zeta
function.

\end{abstract}


\titlerunning{Geometric Polynomials}

\institute{Khristo N. Boyadzhiev \at
              Department of Mathematics and Statistics, Ohio Northern University Ada, Ohio 45810, USA\\
              \email{k-boyadzhiev@onu.edu},
              \and
              Ayhan Dil \at
              Department of Mathematics, Akdeniz University, 07058-Antalya, Turkey\\
              \email{adil@akdeniz.edu.tr}
              }

\section{Introduction}

Let $S(n,k)$ be the Stirling numbers of the second kind (see \cite{B1,B6,GKP}%
). The geometric polynomials
\begin{equation*}
\omega _{n}(x)=\sum_{k=0}^{n}S(n,k)k!x^{k}
\end{equation*}%
were discussed and used in \cite{B2,B4,B5,A1,A2,A3}. These polynomials are
related to the geometric series in the following way
\begin{equation*}
\left( x\frac{d}{dx}\right) ^{m}\frac{1}{1-x}=\sum_{k=0}^{\infty }k^{m}x^{k}=%
\frac{1}{1-x}\omega _{m}\left( \frac{x}{1-x}\right)
\end{equation*}%
for every $|x|\;<1$ and every $m=0,1,2,...$ . Here are the first five of
them
\begin{equation*}
\begin{array}{l}
{\omega _{0}(x)=1} \\
{\omega _{1}(x)=x} \\
{\omega _{2}(x)=2x^{2}+x} \\
{\omega _{3}(x)=6x^{3}+6x^{2}+x} \\
{\omega _{4}(x)=24x^{4}+36x^{3}+14x^{2}+x\,\;.}%
\end{array}%
\end{equation*}%
The polynomials $\omega _{n}\left( x\right) $ can be extended to a more
general form depending on a parameter
\begin{equation}
\omega _{n,r}(x)=\frac{1}{\Gamma (r)}\sum_{k=0}^{n}S(n,k)\Gamma (k+r)x^{k}
\label{0}
\end{equation}%
for every $r>0$, where $\omega _{n,1}\left( x\right) =\omega _{n}\left(
x\right) $ . The polynomials $\omega _{n,r}(x)$ have the property
\begin{align}
\left( x\frac{d}{dx}\right) ^{m}\frac{1}{(1-x)^{r+1}}& =\sum_{k=0}^{\infty }%
\binom{k+r}{k}k^{m}x^{k}  \notag \\
& =\frac{1}{(1-x)^{r+1}}\omega _{m,r+1}\left( \frac{x}{1-x}\right)  \label{1}
\end{align}%
for any $m,r=0,1,2,...$ and also participate in the series transformation
formula
\begin{equation*}
\sum_{k=0}^{\infty }\binom{k+r}{k}f(k)\,x^{k}=\frac{1}{(1-x)^{r+1}}%
\sum_{m=0}^{\infty }\frac{f^{(m)}(0)}{m!}\omega _{m,r+1}\left( \frac{x}{1-x}%
\right)
\end{equation*}%
for appropriate entire functions $f(z)$, see \cite{B2}. The first five
polynomials for $n=0,1,2,3,4$ are:%
\begin{equation*}
\begin{array}{l}
{\omega _{0,r}(x)=1} \\
{\omega _{1,r}(x)=rx} \\
{\omega _{2,r}(x)=r(r+1)x^{2}+rx} \\
{\omega _{3,r}(x)=r(r+1)(r+2)x^{3}+3r(r+1)x^{2}+rx} \\
{\omega
_{4,r}(x)=r(r+1)(r+2)(r+3)x^{4}+6r(r+1)(r+2)x^{3}+7r(r+1)x^{2}+rx\,\;.}%
\end{array}%
\end{equation*}%
We can write $\omega _{n,r}(x)$ also in the form%
\begin{equation}
\omega _{n,r}(x)=\sum_{k=0}^{n}S(n,k)r(r+1)...(r+k-1)x^{k}  \label{2}
\end{equation}%
and from here we find%
\begin{equation*}
\omega _{n,r}(-1)=\sum_{k=0}^{n}S(n,k)(-r)(-r-1)...(-r-k+1)=(-r)^{n}
\end{equation*}%
and%
\begin{equation*}
\omega _{n,r}(1):=\omega _{n,r}=\sum_{k=0}^{n}S(n,k)r(r+1)...(r+k-1).
\end{equation*}

On the other hand, considering the Pochhammer symbol%
\begin{equation*}
\left( x\right) _{n}=x\left( x+1\right) \ldots\left( x+n-1\right) =\frac{%
\Gamma\left( x+n\right) }{\Gamma\left( x\right) }
\end{equation*}
we have%
\begin{equation}
w_{n,r}\left( x\right) =\sum_{k=0}^{n}S\left( n,k\right) \left( r\right)
_{k}x^{k}.  \label{a1}
\end{equation}
Also we know that the Pochhammer symbol can be written in terms of Stirling
numbers of the first kind $s\left( n,k\right) $, as%
\begin{equation}
\left( r\right) _{k}=\sum_{i=0}^{k}\left( -1\right) ^{k-i}s\left( k,i\right)
r^{i}.  \label{a2}
\end{equation}
Hence we can write geometric polynomials in terms of the Stirling numbers of
the first and second kind as:%
\begin{equation}
w_{n,r}\left( x\right) =\sum_{k=0}^{n}\sum_{i=0}^{k}S\left( n,k\right)
s\left( k,i\right) \left( -r\right) ^{i}\left( -x\right) ^{k}.  \label{a3}
\end{equation}
\bigskip

The polynomials $\omega_{n,r}(x)$ can naturally be extended for $r=0$ (as $%
\omega_{n,0}(x)=\delta_{n,0}$) and for $r<0$ by formula (\ref{2}). Thus $%
\omega_{n,-1}(x)=-x$ for all $n\geq1$ . The polynomials $\omega_{n,-r}(x)$
will be formally addressed in Proposition \ref{pr1} below. In what follows, $%
\omega_{n,r}$ are simply called geometric polynomials.

The purpose of this article is to present further properties and
applications of the polynomials $\omega_{n,r}(x)$. To prove some of these
properties we shall use the close relationship of $\omega_{n,r}(x)$ to the
exponential polynomials%
\begin{equation*}
\varphi_{n}(x)=\sum_{k=0}^{n}S(n,k)x^{k}
\end{equation*}
which were studied in \cite{B1,B2,A1}. The first five exponential
polynomials are%
\begin{equation*}
\begin{array}{l}
{\varphi_{0}(x)=1} \\
{\varphi_{1}(x)=x} \\
{\varphi_{2}(x)=x^{2}+x} \\
{\varphi_{3}(x)=x^{3}+3x^{2}+x} \\
{\varphi_{4}(x)=x^{4}+6x^{3}+7x^{2}+x\;.}%
\end{array}%
\end{equation*}

The geometric polynomials originate from the works of Leonhard Euler. They
are proven to be an effective tool in different topics in combinatorics and
analysis. The generalized geometric polynomials help to solve a wider range
of problems, as demonstrated in the present paper and in some previous works
by the authors (see \cite{B1,B2,B4,B6,A2,A3}).

Now a brief summary of the other sections. In the second section we obtain
generating functions for the geometric polynomials $\omega_{n,r}(x)$
according to the parameters $n$ and $r$. We also include there a technical
result involving Lah numbers. Section three contains several recurrence
relations and differential equations for our polynomials. In section four we
list several integral representations of $\omega_{n,r}(x)$ obtained by using
classical formulas. In section five we evaluate in closed form several power
series where the coefficients include values of the Riemann zeta function.
For example, the following series is evaluated in closed form (for any $|x|<2
$ and for arbitrary integers $r\geq0,\;p>0$ )%
\begin{equation*}
\sum_{n=1}^{\infty}\binom{n+r}{n}n^{p}\left\{ \zeta(n+1)-1\right\} x^{n}.
\end{equation*}
Furthermore, we extend these results to series with values of Euler's eta
function and the Lerch Transcendent.

\section{Generating functions}

\begin{lem}
For every $n=0,1,2,...$ and every $r>0$ we have the integral representation
\begin{equation}
\omega_{n,r}(x)=\frac{1}{\Gamma(r)}\int_{0}^{\infty}\lambda^{r-1}\varphi
_{n}(x\lambda)e^{-\lambda}d\lambda.  \label{g0}
\end{equation}
\end{lem}

\begin{prf*}

Evaluating this integral we immediately find%
\begin{align*}
\int_{0}^{\infty}\lambda^{r-1}\varphi_{n}(x\lambda)e^{-\lambda}d\lambda &
=\sum_{k=0}^{n}S(n,k)x^{k}\int_{0}^{\infty}\lambda^{k+\lambda-1}e^{-\lambda
}d\lambda \\
& =\sum_{k=0}^{n}S(n,k)\Gamma(k+r)x^{k}=\Gamma(r)\omega_{n,r}(x).
\end{align*}
\end{prf*}

\begin{prop}
The exponential generating function for $\omega_{n,r}(x)$ is%
\begin{equation}
\left( 1-x(e^{t}-1)\right) ^{-r}=\sum_{n=0}^{\infty}\omega_{n,r}(x)\frac{%
t^{n}}{n!}.  \label{g1}
\end{equation}
In particular, when $r=1$,%
\begin{equation*}
\frac{1}{1-x(e^{t}-1)}=\sum_{n=0}^{\infty}\omega_{n}(x)\frac{t^{n}}{n!}.
\end{equation*}
\end{prop}

\begin{prf*}
Let us consider the well-known generating function for the exponential
polynomials%
\begin{equation*}
e^{x(e^{t}-1)}=\sum_{n=0}^{\infty}\varphi_{n}(x)\frac{t^{n}}{n!}
\end{equation*}
\cite{B1,B2}. From here and Lemma 1,%
\begin{equation*}
\int_{0}^{\infty}\lambda^{r-1}e^{x\lambda(e^{t}-1)}e^{-\lambda}d\lambda
=\Gamma(r)\sum_{n=0}^{\infty}\omega_{n,r}(x)\frac{t^{n}}{n!}.
\end{equation*}
At the same time,%
\begin{equation*}
\int_{0}^{\infty}\lambda^{r-1}e^{x\lambda(e^{t}-1)}e^{-\lambda}d\lambda
=\int_{0}^{\infty}\lambda^{r-1}e^{-\lambda\lbrack1-x(e^{t}-1)]}d\lambda =%
\frac{\Gamma(r)}{[1-x(e^{t}-1)]^{r}}
\end{equation*}
and the proof is completed.
\end{prf*}

\subsection{Generating functions of the $w_{n,r}\left( x\right) $
polynomials with respect to the variable $r$}

The first proposition gives a form of the ordinary generating function of $%
w_{n,r}\left( x\right) $.

\begin{prop}
\begin{equation*}
\sum_{r=0}^{\infty}w_{n,r}\left( x\right) t^{r}=\frac{t}{1-t}w_{n}\left(
\frac{x}{1-t}\right) .
\end{equation*}
\end{prop}

\begin{prf*}
Multiplying the both sides of (\ref{a1}) by\ $t^{r}$ and summing from $r=0$
to $\infty$ we get%
\begin{equation*}
\sum_{r=0}^{\infty}w_{n,r}\left( x\right) t^{r}=\sum_{k=0}^{n}S\left(
n,k\right) x^{k}\sum_{r=0}^{\infty}\left( r\right) _{k}t^{r}.
\end{equation*}
Here using the equation%
\begin{equation*}
\sum_{r=0}^{\infty}\left( r\right) _{k}t^{r}=k!\frac{t}{\left( 1-t\right)
^{k+1}},
\end{equation*}
(see \cite{G,AP}) we have%
\begin{equation*}
\sum_{r=0}^{\infty}w_{n,r}\left( x\right) t^{r}=\frac{t}{1-t}%
\sum_{k=0}^{n}S\left( n,k\right) k!\left( \frac{x}{1-t}\right) ^{k}
\end{equation*}
which completes the proof.
\end{prf*}

We have also a form of the exponential generating function of the $%
w_{n,r}\left( x\right) $ as follows:

\begin{prop}
\label{k1}%
\begin{equation*}
\sum_{r=0}^{\infty}w_{n,r}\left( x\right) \frac{t^{r}}{r!}=e^{t}\left[
\varphi_{n}\left( xt\right) +t\sum_{k=0}^{n}S\left( n,k\right)
P_{k-1}^{k}\left( t\right) x^{k}\right] ,
\end{equation*}
where for $n,k\in\mathbb{N}$%
\begin{equation*}
P_{k}^{n}\left( x\right) =\sum_{j=0}^{k-1}\sum_{m=0}^{k-j-1}\left(
n-m\right) \binom{k}{j}x^{j}.
\end{equation*}
\end{prop}

\begin{prf*}
Considering the equation \cite{AP}%
\begin{equation*}
\sum_{r=0}^{\infty}\left( r\right) _{k}\frac{t^{r}}{r!}=te^{t}\left(
t^{k-1}+P_{k-1}^{k}\left( t\right) \right)
\end{equation*}
and (\ref{a1}),\ we can write%
\begin{equation*}
\sum_{r=0}^{\infty}w_{n,r}\left( x\right) \frac{t^{r}}{r!}%
=\sum_{k=0}^{n}S\left( n,k\right) x^{k}te^{t}\left(
t^{k-1}+P_{k-1}^{k}\left( t\right) \right) .
\end{equation*}
Remembering the definition of exponential polynomials $\varphi_{n}\left( x\right) $  we have%
\begin{equation*}
\sum_{r=0}^{\infty}w_{n,r}\left( x\right) \frac{t^{r}}{r!}=e^{t}\left[
\sum_{k=0}^{n}S\left( n,k\right) \left( xt\right)
^{k}+t\sum_{k=0}^{n}S\left( n,k\right) x^{k}P_{k-1}^{k}\left( t\right) %
\right] ,
\end{equation*}
which is the desired result.
\end{prf*}

Lastly we give a form of the Dirichlet generating function of $w_{n,r}\left(
x\right) $, which is an immediate result from the equation (\ref{a3}).

\begin{prop}
For any real number $\sigma$ such that $\sigma>k+1$ we have%
\begin{equation*}
\sum_{r=1}^{\infty}\frac{w_{n,r}\left( x\right) }{r^{\sigma}}=\sum_{k=0}^{n}%
\left[ \sum_{j=1}^{k}\left( -1\right) ^{j+k}S\left( n,k\right) s\left(
k,j\right) \zeta\left( \sigma-j\right) \right] x^{k}.
\end{equation*}
\end{prop}

Now we give the exponential generating function of the Pochhammer symbol in
terms of the well-known (unsigned) Lah numbers (see sequence A105278 in
OEIS)\ which are defined by%
\begin{equation*}
L(n,k)=\frac{n!}{k!}\binom{n-1}{k-1}.
\end{equation*}

\begin{prop}
We have the following series representation:%
\begin{equation*}
\sum_{r=0}^{\infty}\left( r\right) _{k}\frac{t^{r}}{r!}=e^{t}\sum_{j=0}^{k}L%
\left( k,j\right) t^{j}.
\end{equation*}
\end{prop}

\begin{prf*}
From equation (\ref{a2}) we write%
\begin{equation*}
\left( r\right) _{k}=(-1)^{k}\sum_{i=0}^{k}s(k,i)(-1)^{i}r^{i}.
\end{equation*}
Now summing on $r$ and changing order of summation on the RHS we find%
\begin{equation*}
\sum_{r=0}^{\infty}\left( r\right) _{k}\frac{t^{r}}{r!}=(-1)^{k}\sum
_{i=0}^{k}s(k,i)(-1)^{i}\left\{ \sum_{r=0}^{\infty}r^{i}\frac{t^{r}}{r!}%
\right\} .
\end{equation*}
Then write%
\begin{equation*}
\sum_{r=0}^{\infty}r^{i}\frac{t^{r}}{r!}=e^{t}\sum_{j=0}^{i}S(i,j)t^{j}
\end{equation*}
this is equation (2.9) in \cite{B1} (or equation (2.4) in \cite{B2}). The equation becomes%
\begin{equation*}
\sum_{r=0}^{\infty}\left( r\right) _{k}\frac{t^{r}}{r!}=(-1)^{k}e^{t}%
\sum_{j=0}^{k}t^{j}\left\{ \sum_{i=o}^{k}s(k,i)S(i,j)(-1)^{i}\right\} ,
\end{equation*}
which can equally be written as in the statement by the help of the relation%
\begin{equation*}
L(k,j)=(-1)^{k}\sum_{i=0}^{k}s(k,i)S(i,j)(-1)^{i}
\end{equation*}
(see page 156 in \cite{C}).
\end{prf*}

As an immediate result of (\ref{a1}) we can state the exponential generating
function of the generalized geometric polynomials in terms of Lah numbers.

\begin{cor}
\label{k2}%
\begin{equation*}
\sum_{r=0}^{\infty}w_{n,r}\left( x\right) \frac{t^{r}}{r!}=e^{t}\sum
_{k=0}^{n}\left( \sum_{j=0}^{k}S\left( n,k\right) L\left( k,j\right)
t^{j}\right) x^{k}.
\end{equation*}
\end{cor}

\begin{cor}
We have the following equation for the partial sums of the power series of
Lah numbers,%
\begin{equation*}
\sum_{j=1}^{k-1}L\left( k,j\right) t^{j-1}=P_{k-1}^{k}\left( t\right) .
\end{equation*}
\end{cor}

\begin{prf*}
Comparing RHS of Proposition \ref{k1}\ and Corollary \ref{k2} we see that%
\begin{equation*}
\varphi_{n}\left( xt\right) +t\sum_{k=0}^{n}S\left( n,k\right)
P_{k-1}^{k}\left( t\right) x^{k}=\sum_{k=0}^{n}\left( \sum_{j=0}^{k}S\left(
n,k\right) L\left( k,j\right) t^{j}\right) x^{k}.
\end{equation*}
After some rearrangement we get%
\begin{equation*}
\varphi_{n}\left( xt\right) =\sum_{k=0}^{n}\left( \sum_{j=0}^{k}L\left(
k,j\right) t^{j}-tP_{k-1}^{k}\left( t\right) \right) S\left( n,k\right)
x^{k}.
\end{equation*}
Now taking into account the definition of exponential polynomials and
comparing the coefficients of $x^{k}$ we get the result.
\end{prf*}

\section{Recurrence relations}

Now we give a recurrence relation for the polynomials $\omega_{n,r}\left(
x\right) $ with respect to the variable $n$.

\begin{prop}
We have%
\begin{equation*}
\frac{w_{n+1,r}\left( x\right) +rw_{n,r}\left( x\right) }{r\left( x+1\right)
}=\sum_{k=0}^{n}\binom{n}{k}w_{k,r}\left( x\right) w_{n-k}\left( x\right) .
\end{equation*}
When $x=1$ this becomes%
\begin{equation*}
w_{n+1,r}=\sum_{k=0}^{n-1}\binom{n}{k}2rw_{k,r}w_{n-k}+rw_{n,r}.
\end{equation*}
\end{prop}

\begin{prf*}
Considering the derivative of the generating function of the $w_{n,r}\left(
x\right) $ polynomials (\ref{g1})\ we have%
\begin{equation*}
\frac{rxe^{t}}{\left( 1-x\left( e^{t}-1\right) \right) ^{r+1}}=\sum
_{n=0}^{\infty}w_{n+1,r}\left( x\right) \frac{t^{n}}{n!}.
\end{equation*}
From this we get%
\begin{equation*}
r\left( \frac{1+x}{\left( 1-x\left( e^{t}-1\right) \right) }-1\right) \frac{1%
}{\left( 1-x\left( e^{t}-1\right) \right) ^{r}}=\sum_{n=0}^{\infty
}w_{n+1,r}\left( x\right) \frac{t^{n}}{n!}.
\end{equation*}
Expanding the LHS and comparing coefficients of both sides completes the proof.
\end{prf*}

\begin{prop}
For any two positive integers $n$ and $m$ we have,%
\begin{equation*}
\omega_{n+m,r}(x)=\sum_{k=0}^{n}\sum_{j=0}^{m}\binom{n}{k}S\left( m,j\right)
\left( r\right) _{j}j^{n-k}x^{j}\omega_{k,r+j}(x).
\end{equation*}
\end{prop}

\begin{prf*}
From (\ref{g0})\ we have%
\begin{equation*}
\omega_{n+m,r}(x)=\frac{1}{\Gamma(r)}\int_{0}^{\infty}\lambda^{r-1}%
\varphi_{n+m}(x\lambda)e^{-\lambda}d\lambda.
\end{equation*}
In the light of the equation (see \cite{B1})%
\begin{equation*}
\varphi_{n+m}\left( x\right) =\sum_{k=0}^{n}\sum_{j=0}^{m}\binom{n}{k}%
S\left( m,j\right) j^{n-k}x^{j}\varphi_{k}\left( x\right)
\end{equation*}
we can calculate the integral on the RHS which completes the proof.
\end{prf*}

\begin{prop}
\label{p3}For every $r>0$ and every $n=0,1,2,...$ , we have the differential
equation%
\begin{equation*}
\omega_{n,r+1}(x)=\frac{x}{r}\omega_{n,r}^{\prime}(x)+\omega_{n,r}(x).
\end{equation*}
In particular, for $r=1$,%
\begin{equation*}
\omega_{n,2}(x)=x\omega_{n}^{\prime}(x)+\omega_{n}(x).
\end{equation*}
\end{prop}

\begin{prf*}
\begin{align*}
\Gamma(r+1)\omega_{n,r+1}(x) & =\sum_{k=0}^{n}S(n,k)\Gamma(k+r+1)x^{k} \\
& =\sum_{k=0}^{n}S(n,k)\Gamma(k+r)(k+r)x^{k} \\
& =\sum_{k=0}^{n}S(n,k)\Gamma(k+r)kx^{k}+r\sum_{k=0}^{n}S(n,k)\Gamma
(k+r)x^{k} \\
& =\Gamma(r)x\omega_{n,r}^{\prime}(x)+r\Gamma(r)\omega_{n,r}(x).
\end{align*}
Dividing by $\Gamma(r+1)$ we come to the desired equation.
\end{prf*}

\begin{rem}
Using the two equations in the above proposition we find immediately%
\begin{equation*}
\omega_{n,3}(x)=\frac{x^{2}}{2}\omega_{n}^{\prime\prime}(x)+2x\omega
_{n}^{\prime}(x)+\omega_{n}(x)
\end{equation*}
and this process can be continued further. Therefore, every polynomial $%
\omega_{n,r}(x)$, where $r>0$ is an integer, can be written in terms of $%
\omega_{n}(x)$ and its derivatives with easily computable special
coefficients.
\end{rem}

\begin{prop}
For every $r>0$ and every $n=0,1,2,...$ , we have the recurrences%
\begin{align*}
\omega_{n+1,r}(x) & =r\left( (x+1)\omega_{n,r+1}(x)-\omega_{n,r}(x)\right) ,
\\
\omega_{n+1,r}(x) & =(x^{2}+x)\omega_{n,r}^{\prime}(x)+xr\omega_{n,r}(x)
\end{align*}
and in particular, for $r=1$,%
\begin{equation*}
\omega_{n+1}(x)=(x^{2}+x)\omega_{n}^{\prime}(x)+x\omega_{n}(x).
\end{equation*}
\end{prop}

\begin{prf*}
Using the following property of exponential polynomials \cite{B1}%
\begin{equation*}
\varphi_{n+1}(x)=x\left( \varphi_{n}(x)+\varphi_{n}^{\prime}(x)\right)
\end{equation*}
we have
\begin{align*}
\Gamma(r)\omega_{n+1,r}(x) & =\int_{0}^{\infty}\lambda^{r-1}\varphi
_{n+1}(x\lambda)\,e^{-\lambda}d\lambda \\
&
=x\int_{0}^{\infty}\lambda^{r}\varphi_{n}(x\lambda)e^{-\lambda}d\lambda+%
\int_{0}^{\infty}\lambda^{r}\varphi_{n}^{\prime}(x\lambda
)e^{-\lambda}d\lambda \\
& =x\Gamma(r+1)\omega_{n,r+1}(x)-\int_{0}^{\infty}\varphi_{n}(x\lambda
)(r\lambda^{r-1}e^{-\lambda}-\lambda^{r}e^{-\lambda})d\lambda
\end{align*}
and this becomes%
\begin{equation*}
\Gamma(r)\omega_{n+1,r}(x)=x\Gamma(r+1)\omega_{n,r+1}(x)-r\Gamma
(r)\omega_{n,r}(x)+\Gamma(r+1)\omega_{n,r+1}(x).
\end{equation*}
Dividing both sides in this equation by $\Gamma(r)$ yields the first
equation in the proposition. Applying Proposition \ref{p3} to $%
\omega_{n,r+1} $ on the RHS brings to the second equation.
\end{prf*}

\begin{rem}
It is interesting that the second equation can also be written in the form%
\begin{equation*}
\omega_{n+1,r}(x)=\varphi_{2}(x)\omega_{n,r}^{\prime}(x)+\varphi_{1}(x)r%
\omega_{n,r}(x).
\end{equation*}
\end{rem}

\begin{prop}
For every $r\geq1$ and $n=0,1,2,...$ , the binomial transform of the
geometric polynomials is given by
\begin{equation*}
\sum_{k=0}^{n}\binom{n}{k}\omega_{k,r}(x)=\left( 1+\frac{1}{x}\right)
\omega_{n,r}(x)-\frac{1}{x}\omega_{n,r-1}(x),
\end{equation*}
where for $r=1$,%
\begin{equation*}
\sum_{k=0}^{n}\binom{n}{k}\omega_{k}(x)=\left( 1+\frac{1}{x}\right)
\omega_{n}(x)-\frac{1}{x}\delta_{n,0},
\end{equation*}
and for $n>0$ this is simply%
\begin{equation*}
\sum_{k=0}^{n}\binom{n}{k}\omega_{k}(x)=\left( 1+\frac{1}{x}\right)
\omega_{n}(x).
\end{equation*}
\end{prop}

\begin{prf*}
For the proof we use a formula from \cite{B1},%
\begin{equation*}
\sum_{k=0}^{n}\binom{n}{k}\varphi_{k}(x)=\varphi_{n}(x)+\varphi_{n}^{\prime
}(x).
\end{equation*}
Using integration by parts we write%
\begin{align*}
& \Gamma(r)\sum_{k=0}^{n}\binom{n}{k}\omega_{k,\,r}(x) \\
& \quad=\int_{0}^{\infty}\lambda^{r-1}e^{-\lambda}\left\{ \varphi
_{n}(x\lambda)+\varphi_{n}^{\prime}(x\lambda)\right\} d\lambda \\
& \quad=\Gamma(r)\omega_{n,\,r}(x)+\frac{1}{x}\int_{0}^{\infty}\lambda
^{r-1}e^{-\lambda}d\,\varphi_{n}^{\prime}(x\lambda) \\
& \quad=\Gamma(r)\omega_{n,\,r}(x)+\frac{1}{x}\left\{ \left. \lambda
^{r-1}e^{-\lambda}\varphi_{n}(x\lambda)\right\vert _{0}^{\infty}-\int
_{0}^{\infty}\varphi_{n}(x\lambda)\left\{ (r-1)\lambda^{r-2}e^{-\lambda
}-\lambda^{r-1}e^{-\lambda}\right\} d\lambda\right\} \\
& \quad=\Gamma(r)\omega_{n,\,r}(x)+\frac{1}{x}\left\{ -(r-1)\Gamma
(r-1)\omega_{n,\,r-1}(x)+\Gamma(r)\omega_{n,\,r}(x)\right\} \\
& \quad=\Gamma(r)\left\{ \omega_{n,\,r}(x)-\frac{1}{x}\omega_{n,\,r-1}(x)+%
\frac{1}{x}\omega_{n,\,r}(x)\right\} ,
\end{align*}
which is the desired result.

When $r=1$ the computation is simplified to%
\begin{align*}
\sum_{k=0}^{n}\binom{n}{k}\omega_{n}(x) & =\int_{0}^{\infty}e^{-\lambda
}\left\{ \varphi_{n}(x\lambda)+\varphi_{n}^{\prime}(x\lambda)\right\}
d\lambda \\
& =\omega_{n}(x)+\frac{1}{x}\left\{ \left.
e^{-\lambda}\varphi_{n}(x\lambda)\right\vert
_{0}^{\infty}+\int_{0}^{\infty}e^{-\lambda}\varphi
_{n}(x\lambda)\,d\lambda\right\} \\
& =\omega_{n}(x)+\frac{1}{x}\left\{ -\varphi_{n}(0)+\omega_{n}(x)\right\} ,
\end{align*}
and this explains the term $\delta_{n,0}$, as $\varphi_{n}(0)=0$ for $n>1$
and $\varphi_{0}(0)=1$.
\end{prf*}

\section{Integral representations involving the geometric polynomials}

In the next proposition we give several integral representations involving
the geometric polynomials. The first one provides a Mellin integral
representation in the general case. In the second representation we use the
Riemann zeta function $\zeta(s)$. For the following proposition we shall use
the well-known estimate for the Gamma function:
\begin{equation*}
\left\vert \Gamma \left( x+iy\right) \right\vert \sim \sqrt{2\pi }\left\vert
y\right\vert ^{x-\frac{1}{2}}e^{-x-\frac{\pi }{2}\left\vert y\right\vert }
\end{equation*}
($\left\vert y\right\vert \rightarrow \infty $) for any fixed real $x$. This
explains the behavior of the Gamma function on vertical lines.

\begin{prop}
For every $r\geq0$, every $0<x<1$, and every $n=0,1,2,...$ we have%
\begin{equation*}
\frac{1}{(1+x)^{r+1}}\omega_{n,r+1}\left( \frac{-x}{1+x}\right) =\frac{%
(-1)^{n}}{2\pi i\Gamma(r+1)}\int_{a-i\infty}^{a+i\infty}x^{-s}s^{n}\Gamma(s)%
\Gamma(r+1-s)ds.
\end{equation*}
Integration here is on a vertical line $\{s=a+it,\;-\infty<t<+\infty\}$,
where $0<a<1$.

For all $x>0$, $n=0,1,...$ and every $a>n+1$,%
\begin{equation*}
\frac{e^{x}}{e^{x}-1}\omega_{n}\left( \frac{1}{e^{x}-1}\right) =\frac {1}{%
2\pi i}\int_{a-i\infty}^{a+i\infty}x^{-s}\zeta(s-n)\Gamma(s)ds.
\end{equation*}
Also, for all $x>0$ and $n=1,2,3,...$%
\begin{align*}
\frac{e^{x}}{e^{x}+1}\omega_{n}\left( \frac{-1}{e^{x}+1}\right) &
=(-1)^{n}\int_{0}^{\infty}\sin\left( xt+\frac{\pi n}{2}\right) \frac{t^{n}}{%
\sinh(\pi t)}dt, \\
\frac{e^{x}}{e^{x}-1}\omega_{n}\left( \frac{1}{e^{x}-1}\right) & =\frac{n!}{%
x^{n+1}}+2(-1)^{n}\int_{0}^{\infty}\sin\left( xt+\frac{\pi n}{2}\right)
\frac{t^{n}}{e^{2\pi t}-1}dt.
\end{align*}
\end{prop}

(For the last representation cf. Ramanujan's Entry 2 on p.335 in his
notebooks; \cite[p. 411]{BB}).

\begin{prf*}
Starting from the Mellin integral representation (formula 5.37 on p.196 in
\cite{F})%
\begin{equation*}
\frac{1}{(1+x)^{r+1}}=\frac{1}{2\pi i\Gamma(r+1)}\int_{a-i\infty}^{a+i\infty
}x^{-s}\Gamma(s)\Gamma(r+1-s)ds,
\end{equation*}
we compute%
\begin{equation*}
\left( x\frac{d}{dx}\right) ^{m}\frac{1}{(1+x)^{r+1}}=\frac{(-1)^{m}}{2\pi
i\Gamma(r+1)}\int_{a-i\infty}^{a+i\infty}x^{-s}s^{m}\Gamma(s)\Gamma(r+1-s)ds
\end{equation*}
and at the same time,%
\begin{align*}
\left( x\frac{d}{dx}\right) ^{m}\frac{1}{(1+x)^{r+1}} & =\left( x\frac {d}{dx%
}\right) ^{m}\sum_{n=0}^{\infty}\binom{r+n}{n}(-x)^{n} \\
& =\sum_{n=0}^{\infty}\binom{r+n}{n}n^{m}(-x)^{n} \\
& =\frac{1}{(1+x)^{r+1}}\omega_{m,r+1}\left( \frac{-x}{1+x}\right) .
\end{align*}
Equating the right hand sides completes the proof of the first
representation. For the second representation we start from the well-known
formula%
\begin{equation*}
e^{-x}=\frac{1}{2\pi i}\int_{a-i\infty}^{a+i\infty}x^{-s}\Gamma(s)ds,
\end{equation*}
where we replace $x$ by $xk$ and multiply both sides by $k^{n}$ in order to
get%
\begin{equation*}
k^{n}e^{-xk}=\frac{1}{2\pi i}\int_{a-i\infty}^{a+i\infty}x^{-s}\frac {1}{%
k^{\,s-n}}\Gamma(s)ds.
\end{equation*}
Summing for $k=1,2,...$ now yields the desired representation. The last two
representations follow from the two integral formulas (both are sine Fourier
transforms)%
\begin{align*}
& \sum_{k=0}^{\infty}(-1)^{k}e^{-kx}=\frac{1}{1+e^{-x}}=\frac{1}{2}+\int
_{0}^{\infty}\frac{\sin(xt)}{\sinh(\pi t)}dt, \\
& \sum_{k=0}^{\infty}e^{-kx}=\frac{1}{1-e^{-x}}=\frac{1}{x}+\frac{1}{2}%
+2\int_{0}^{\infty}\frac{\sin(xt)}{e^{2\pi t}-1}dt.
\end{align*}
Differentiation $n$ times for $x$ yields the desires representations.
\end{prf*}

Setting $x\rightarrow0$ in the last two representations we obtain the
corollary.

\begin{cor}
For every $n=1,2,3,...$ we have%
\begin{align*}
\omega_{n}\left( \frac{-1}{2}\right) & =2(-1)^{n}\sin\frac{\pi n}{2}%
\int_{0}^{\infty}\frac{t^{n}}{\sinh(\pi t)}dt \\
& =\frac{4(-1)^{n}n!}{\pi^{n+1}}\left( 1-\frac{1}{2^{n+1}}\right) \sin \frac{%
\pi n}{2}\zeta(n+1)
\end{align*}
and%
\begin{align*}
& \mathop{\lim }\limits_{x\,\rightarrow0}\left\{ \frac{e^{x}}{e^{x}-1}%
\omega_{n}\left( \frac{1}{e^{x}-1}\right) -\frac{n!}{x^{n+1}}\right\} \\
& =2(-1)^{n}\sin\frac{\pi n}{2}\int_{0}^{\infty}\frac{t^{n}}{e^{2\pi t}-1}dt=%
\frac{2n!(-1)^{n}}{(2\pi)^{n+1}}\sin\frac{\pi n}{2}\zeta(n+1).
\end{align*}
In particular, $\omega_{2k}\left( \frac{-1}{2}\right) =0,\;\;k=1,2,...\;.$
\end{cor}

\begin{prop}
\label{pr1}Let $r>0$. Defining for every $n\geq0$ the polynomials%
\begin{equation*}
\omega_{n,-r}(x)=\sum_{k=0}^{n}S(n,k)(-1)^{k}r(r-1)...(r-k+1)x^{k},
\end{equation*}
we have for $p=0,1,2,...$
\begin{equation*}
\left( x\frac{d}{dx}\right) ^{p}(1-x)^{r}=(1-x)^{r}\omega_{p,-r}\left( \frac{%
x}{1-x}\right)
\end{equation*}
and when $|x|<1$ this can also be written as%
\begin{equation}
\sum_{n=0}^{\infty}\binom{r}{n}n^{p}x^{n}=(1+\,x)^{r}\,%
\sum_{k=0}^{p}S(p,k)r(r-1)...(r-k+1)\left( \frac{x}{1+\,x}\right) ^{k}.
\label{3}
\end{equation}
(Changing $x$ to $-x,$ expanding $(1+x)^{r}$ by the binomial formula, and
applying $\left( x\frac{d}{dx}\right) ^{p}$\ to both sides). When $r$ is not
an integer, (\ref{3}) represents a closed form evaluation of the infinite
series on the LHS.
\end{prop}

\begin{prf*}
We use the well-known formula%
\begin{equation*}
\left( x\frac{d}{dx}\right) ^{p}f(x)=\sum_{k=0}^{p}S(p,k)x^{k}f^{(k)}(x)
\end{equation*}
for every $p$-times differentiable function \cite{B1} in order to compute%
\begin{equation*}
\left( x\frac{d}{dx}\right)
^{p}(1-x)^{r}=\sum_{k=0}^{p}S(p,k)x^{k}(-1)^{k}r(r-1)...(r-k+1)(1-x)^{r-k}
\end{equation*}
as needed.
\end{prf*}

\section{\textbf{Applications. Series with zeta values and other series}}

Some applications were given in \cite{B2,B3,A1,A2}. Here we present
additional examples.

\textbf{Example} For $n=0,1,2,...$ considering $s(n,k)$ the Stirling numbers
of the first kind \cite{GKP} with generating polynomial%
\begin{equation*}
f(x)=\binom{x}{p}=\frac{1}{p!}\sum_{k=0}^{p}s(p,k)x^{k}.
\end{equation*}
Applying for this function the transformation formula%
\begin{equation*}
\sum_{k=0}^{\infty}\binom{k+r}{k}f(k)x^{k}=\frac{1}{(1-x)^{r+1}}\sum
_{m=0}^{\infty}\frac{f^{(m)}(0)}{m!}\omega_{m,r+1}\left( \frac{x}{1-x}\right)
\end{equation*}
we find the closed form evaluation%
\begin{align*}
\sum_{k=0}^{\infty}\binom{k+r}{k}\binom{k}{p}x^{k} & =\frac{1}{p!(1-x)^{r+1}}%
\sum_{m=0}^{p}s(p,m)\omega_{m,r+1}\left( \frac{x}{1-x}\right) \\
& =\frac{1}{p!\Gamma(r+1)(1-x)^{r+1}}\sum_{m=0}^{p}s(p,m)%
\sum_{k=0}^{m}S(m,k)\Gamma(k+r+1)\left( \frac{x}{1-x}\right) ^{k} \\
& =\frac{1}{p!\Gamma(r+1)(1-x)^{r+1}}\sum_{k=0}^{p}\Gamma(k+r+1)\left( \frac{%
x}{1-x}\right) ^{k}\sum_{m=0}^{p}s(p,m)S(m,k) \\
& =\frac{\Gamma(p+r+1)}{p!\Gamma(r+1)(1-x)^{r+1}}\left( \frac{x}{1-x}\right)
^{p}=\binom{p+r}{p}\left( \frac{x}{1-x}\right) ^{p}\frac {1}{(1-x)^{r+1}},
\end{align*}
as $\sum_{m=0}^{p}s(p,m)S(m,k)=\delta_{p,k}$ . That is, we obtained the
identity%
\begin{equation*}
\sum_{k=0}^{\infty}\binom{k+r}{k}\binom{k}{p}x^{k}=\binom{p+r}{p}\frac{x^{p}%
}{(1-x)^{p+r+1}}.
\end{equation*}
Once this formula is discovered it can be given a direct proof based on the
well-known expansion \cite{G}%
\begin{equation*}
\sum_{k=0}^{\infty}\binom{k+r}{k}\,x^{k}=\frac{1}{(1-x)^{r+1}}.
\end{equation*}
In the following propositions we extend some results for series with zeta
values including a result of Adamchik and Srivastava (see pp. 142-156 in
\cite{SC}). We shall use the Hurwitz zeta function%
\begin{equation*}
\zeta(s,a)=\sum_{k=0}^{\infty}\frac{1}{(k+a)^{s}},
\end{equation*}
($Re\,\,a>0,\;\,Re\,\,s>1$) and the Riemann zeta function $\zeta
(s)=\zeta(s,1)$.

\begin{prop}
For every $r\geq0$, every integer $p\geq0$ with $r+p>0$, and every $|x|<2$,%
\begin{equation*}
\sum_{n=1}^{\infty}\binom{n+r}{n}n^{p}\left\{ \zeta(n+r+1)-1\right\} x^{n}=%
\frac{1}{\Gamma(r+1)}\sum_{j=0}^{p}S(p,j)\Gamma(r+j+1)\zeta
(r+j+1,2-x)x^{j}\,.
\end{equation*}
When $x=1$ this becomes%
\begin{equation*}
\sum_{n=1}^{\infty}\binom{n+r}{n}n^{p}\left\{ \zeta(n+r+1)-1\right\} =\frac{1%
}{\Gamma(r+1)}\sum_{j=0}^{p}S(p,j)\,\Gamma(r+j+1)\zeta(r+j+1)\,\,.
\end{equation*}
When $r=0$,%
\begin{equation*}
\sum_{n=1}^{\infty}n^{p}\left\{ \zeta(n+1)-1\right\}
x^{n}=\sum_{j=1}^{p}S(p,j)\Gamma(j+1)\zeta(j+1,2-x)x^{j}\,
\end{equation*}
with the summation on the RHS starting from $j=1$, as $p>0$ and $S(p,0)=0$ .
\end{prop}

When $p=0$ and $r>0$ we start the summation on the LHS from $n=0$ to get
\begin{equation*}
\sum_{n=0}^{\infty}\binom{n+r}{n}\left\{ \zeta(n+r+1)-1\right\}
x^{n}=\zeta(r+1,2-x),
\end{equation*}
which is equation (18) on p.146 in \cite{SC}.

In the case $r=p=0$ and $|x|<1$ the series is reduced to the well-known (%
\cite{SC})
\begin{equation*}
\sum_{n=1}^{\infty}\,\zeta(n+1)x^{n}=-\psi(1-x)-\gamma,
\end{equation*}
where $\psi(z)$ is the digamma function and $\gamma$ is Euler's constant.
This case is included in Proposition \ref{p11} below.

\begin{prf*}
We compute%
\begin{align*}
& \sum_{n=1}^{\infty}\binom{n+r}{n}n^{p}\left\{ \zeta(n+r+1)-1\right\} x^{n}
\\
& =\sum_{n=1}^{\infty}\binom{n+r}{n}n^{p}x^{n}\left\{ \sum_{k=1}^{\infty }%
\frac{1}{(k+1)^{n+r+1}}\right\} \\
& =\sum_{k=1}^{\infty}\frac{1}{(k+1)^{r+1}}\left\{ \sum_{n=1}^{\infty}\binom{%
n+r}{n}n^{p}\left( \frac{x}{k+1}\right) ^{n}\right\} \\
& =\sum_{k=1}^{\infty}\frac{1}{(k+1)^{r+1}}\left\{ \left( \frac{k+1}{k+1-x}%
\right) ^{r+1}\omega_{p,r+1}\left( \frac{x}{k+1-x}\right) \right\} \\
& =\sum_{k=1}^{\infty}\frac{1}{(k+1-x)^{r+1}}\omega_{p,r+1}\left( \frac {x}{%
k+1-x}\right)
\end{align*}
by applying formula (\ref{1}) with $\frac{x}{k+1}$ in the place of $x$ .
Further, considering (\ref{0})\ this equals%
\begin{equation*}
\frac{1}{\Gamma(r+1)}\sum_{k=1}^{\infty}\frac{1}{(k+1-x)^{r+1}}\left\{
\sum_{j=0}^{p}S(p,j)\Gamma(r+j+1)\frac{x^{j}}{(k+1-x)^{j}}\right\} \,\,
\end{equation*}%
\begin{equation*}
=\frac{1}{\Gamma(r+1)}\sum_{j=0}^{p}S(p,j)\Gamma(r+j+1)x^{j}\sum_{k=1}^{%
\infty}\frac{1}{(k+1-x)^{r+j+1}}
\end{equation*}%
\begin{equation*}
=\frac{1}{\Gamma(r+1)}\sum_{j=0}^{p}S(p,j)\Gamma(r+j+1)\zeta(r+j+1,2-x)x^{j}%
\,,
\end{equation*}
and the proof is complete.
\end{prf*}

\begin{rem}
In the above proposition if we start with the series%
\begin{equation*}
\sum_{n=1}^{\infty}\binom{n+r}{n}n^{p}\zeta(n+r+1)x^{n}
\end{equation*}
under the condition $|x|<1,$\ then we obtain%
\begin{equation*}
\sum_{n=1}^{\infty}\binom{n+r}{n}n^{p}\zeta(n+r+1)x^{n}=\frac{1}{\Gamma (r+1)%
}\sum_{j=0}^{p}S(p,j)\Gamma(r+j+1)\zeta(r+j+1,1-x)x^{j}.
\end{equation*}
\end{rem}

The next result is based on the same idea.

\begin{prop}
\label{p10}For every integers $r\geq0$, $p>0$, and every $|x|<2$,%
\begin{equation*}
\sum_{n=1}^{\infty}\binom{n+r}{n}n^{p}\left\{ \zeta(n+1)-1\right\} x^{n}=%
\frac{1}{r!}\sum_{j=0}^{p}S(p,j)(r+j)!A(r,x,j)x^{j}\,,
\end{equation*}
where%
\begin{equation*}
A\,(r,x,j)=\sum_{m=0}^{r}\binom{r}{m}x^{m}\zeta(m+j+1,2-x).
\end{equation*}
\end{prop}

\begin{prf*}
Starting as before we find%
\begin{align*}
& \sum_{n=1}^{\infty}\binom{n+r}{n}n^{p}\left\{ \zeta(n+1)-1\right\} x^{n} \\
& =\sum_{n=1}^{\infty}\binom{n+r}{n}n^{p}x^{n}\left\{ \sum_{k=1}^{\infty }%
\frac{1}{(k+1)^{n+1}}\right\} \\
& =\sum_{k=1}^{\infty}\frac{1}{k+1}\left\{ \sum_{n=1}^{\infty}\binom{n+r}{n}%
n^{p}\left( \frac{x}{k+1}\right) ^{n}\right\} \\
& =\sum_{k=1}^{\infty}\frac{1}{k+1}\left\{ \frac{(k+1)^{r+1}}{(k+1-x)^{r+1}}%
\omega_{p,r+1}\left( \frac{x}{k+1}\right) \right\} \\
& =\frac{1}{\Gamma(r+1)}\sum_{j=0}^{p}S(p,j)\Gamma(j+r+1)\,x^{j}\left\{
\sum_{k=1}^{\infty}\frac{(k+1)^{r}}{(k+1-x)^{r+j+1}}\right\} .
\end{align*}
We compute now the sums%
\begin{align*}
A\,(r,x,j) & =\sum_{k=1}^{\infty}\frac{(k+1)^{r}}{(k+1-x)^{r+j+1}}%
=\sum_{k=1}^{\infty}\frac{(k+1-x+x)^{r}}{(k+1-x)^{r+j+1}} \\
& =\sum_{k=1}^{\infty}\frac{1}{(k+1-x)^{r+j+1}}\sum_{m=0}^{r}\binom{r}{m}%
(k+1-x)^{r-m}x^{m} \\
& =\sum_{m=0}^{r}\binom{r}{m}x^{m}\left\{ \sum_{k=1}^{\infty}\frac {1}{%
(k+1-x)^{m+j+1}}\right\} \\
& =\sum_{m=0}^{r}\binom{r}{m}x^{m}\zeta(m+j+1,2-x),
\end{align*}
and the proof is finished.
\end{prf*}

\begin{rem}
If we start with $\zeta(n+1),$ for every integers $r\geq0$, $p>0$, and every
$|x|<1$ we have%
\begin{equation*}
\sum_{n=1}^{\infty}\binom{n+r}{n}n^{p}\zeta(n+1)x^{n}=\frac{1}{r!}\sum
_{j=0}^{p}S(p,j)(r+j)!B(r,x,j)x^{j}\,\,,
\end{equation*}
where%
\begin{equation*}
B(r,x,j)=\sum_{m=0}^{r}\binom{r}{m}x^{m}\zeta(m+j+1,1-x).
\end{equation*}
\end{rem}

For completeness we consider also the case when $p=0$ and $r\geq0$ is an
integer.

\begin{prop}
\label{p11}For every integer $r\geq0$ and every $|x|<1$,%
\begin{equation*}
\sum_{n=1}^{\infty}\binom{n+r}{n}\zeta(n+1)x^{n}=\sum_{m=0}^{r}\binom{r}{m}%
\zeta(m+1,1-x)x^{m}-\psi(1-x)-\gamma.
\end{equation*}
\end{prop}

\begin{prf*}
\begin{align*}
\sum_{n=1}^{\infty}\binom{n+r}{n}\zeta(n+1)x^{n} & =\sum_{n=1}^{\infty }%
\binom{n+r}{n}x^{n}\left\{ \sum_{k=1}^{\infty}\frac{1}{k^{n+1}}\right\} \\
& =\sum_{k=1}^{\infty}\frac{1}{k}\left\{ \sum_{n=1}^{\infty}\binom{n+r}{n}%
\left( \frac{x}{k}\right) ^{n}\right\} \\
& =\sum_{k=1}^{\infty}\frac{1}{k}\left\{ \left( \frac{k}{k-x}\right)
^{r+1}-1\right\} \\
& =\sum_{k=1}^{\infty}\left\{ \frac{(k-x+x)^{r}}{(k-x)^{r+1}}-\frac{1}{k}%
\right\} .
\end{align*}
Writing now
\begin{equation*}
(k-x+x)^{r}=\sum_{m=0}^{r}\binom{r}{m}(k-x)^{r-m}x^{m}
\end{equation*}
we continue the above equation to%
\begin{align*}
& \sum_{k=1}^{\infty}\left\{ \sum_{m=1}^{r}\binom{r}{m}\frac{1}{(k-x)^{m+1}}+%
\frac{1}{k-x}-\frac{1}{k}\right\} \\
& =\sum_{m=1}^{r}\binom{r}{m}\left\{ \sum_{k=1}^{\infty}\frac{1}{(k-x)^{m+1}}%
\right\} +\sum_{k=1}^{\infty}\frac{x}{k(k-x)} \\
& =\sum_{m=0}^{r}\binom{r}{m}\zeta(m+1,1-x)x^{m}-\psi(1-x)-\gamma,
\end{align*}
as (see \cite{SC})%
\begin{equation*}
\psi(1+z)+\gamma=\sum_{k=1}^{\infty}\frac{z}{k(k+z)}
\end{equation*}
and the proof is finished.
\end{prf*}

Similar results hold for the functions%
\begin{equation*}
\eta(s,a)=\sum_{n=0}^{\infty}\frac{(-1)^{n}}{(n+a)^{s}}\text{ and }%
\eta(s)=\eta(s,1),
\end{equation*}
where $Re\left( a\right) >0,\;Re\left( s\right) >0$. The first one is
sometimes called Lerch's eta function and the second one (often used by
Euler) is Euler's eta function.

The next two propositions are proved the same way. We leave details to the
reader.

\begin{prop}
For every $r\geq0$, every integer $p\geq0$, and every $|x|<2$,%
\begin{equation*}
\sum_{n=1}^{\infty}\binom{n+r}{n}n^{p}\left\{ 1-\eta(n+r+1)\right\} x^{n}=%
\frac{1}{\Gamma(r+1)}\sum_{j=0}^{p}S(p,j)\Gamma(r+j+1)\eta (r+j+1,2-x)x^{j}.
\end{equation*}
When $r=p=0$ with summation on the LHS starting from $n=0$ we have
\begin{equation*}
\sum_{n=0}^{\infty}\left\{ 1-\eta(n+1)\right\} x^{n}=\eta(1,2-x).
\end{equation*}
\end{prop}

\begin{rem}
For every $r\geq0$, every integer $p\geq0$, and every $|x|<1$ we have,%
\begin{equation*}
\sum_{n=1}^{\infty}\binom{n+r}{n}n^{p}\eta(n+r+1)x^{n}=\frac{1}{\Gamma (r+1)}%
\sum_{j=0}^{p}S(p,j)\Gamma(r+j+1)\eta(r+j+1,1-x)x^{j},
\end{equation*}
and also%
\begin{equation*}
\sum_{n=1}^{\infty}\binom{n+r}{n}n^{p}\eta(n+1)x^{n}=\frac{1}{\Gamma(r+1)}%
\sum_{j=0}^{p}S(p,j)\Gamma(r+j+1)x^{j}C\left( r,x,j\right)
\end{equation*}
where%
\begin{equation*}
C(r,x,j)=\sum_{m=0}^{r}\binom{r}{m}x^{m}\eta(m+j+1,1-x).
\end{equation*}
\end{rem}

This result can be extended further by involving the Lerch Transcendent%
\begin{equation*}
\Phi(s,z,a)=\sum_{n=0}^{\infty}\frac{z^{n}}{(n+a)^{s}}
\end{equation*}
assuming $a>0,\;\;|z|\;\leq1,\;z\neq1$. Repeating the steps in the proof of
Proposition \ref{p10} one can obtain the following result.

\begin{prop}
For every integer $r\geq0$, every integer $p\geq0$, every $|x|<a$, and every
$|z|\leq1,\;\,z\neq1$ we have%
\begin{equation*}
\sum_{n=0}^{\infty}\binom{n+r}{n}n^{p}\Phi(n+r+1,z,a)x^{n}=\frac{1}{r!}%
\sum_{j=0}^{p}S(p,j)\Gamma(j+r+1)x^{j}\Phi\left( r+j+1,z,a-x\right) .
\end{equation*}
\end{prop}

\begin{rem}
As the previous remarks, for every integer $r\geq 0$, every integer $p\geq 0$%
, every $|x|<a$, and every $|z|\leq 1,\;\,z\neq 1$ we have%
\begin{equation*}
\sum_{n=0}^{\infty }\binom{n+r}{n}n^{p}\Phi (n,z,a)x^{n}=\frac{1}{r!}%
\sum_{j=0}^{p}S(p,j)\Gamma (j+r+1)x^{j}\Phi \left( j,z,a-x\right) .
\end{equation*}
\end{rem}

\begin{acknowledgement}
The authors are grateful to the reviewer for a number of valuable comments.
\end{acknowledgement}

\end{document}